\documentclass{amsart}

\usepackage{amsfonts, amssymb}

\newtheorem{theorem}{Theorem}[section]

\newtheorem{proposition}[theorem]{Proposition}

\theoremstyle{definition}
\newtheorem{definition}[theorem]{Definition}
\theoremstyle{remark}
\newtheorem{remark}[theorem]{Remark}

\numberwithin{equation}{section}

\newcommand{\CC}{\mathbb C}
\newcommand{\RR}{\mathbb R}

\newcommand{\aut}[1]{\text{Aut\,}(#1)}

\newcommand{\supp}{\text{\rm supp }}
\def\cA{{\mathcal A}}

\def\cC{{\mathcal C}}

\def\cH{{\mathcal H}}

\def\im{\textrm{Im}\,}
\def\re{\textrm{Re}\,}

\begin{document}

\title[Invariant metrics of Kohn-Nirenberg domains]
{Bergman and Caratheodory metrics of \\ the 
Kohn-Nirenberg domains}

\author{Taeyong Ahn, Herv\'e Gaussier, 
Kang-Tae Kim}
\address{(Ahn) Center for Geometry and its Applications, 
POSTECH, Pohang 790-784, The Republic of Korea}
\email{(Ahn) triumph@postech.ac.kr}
\address{(Gaussier) Institut Fourier, Universit\'e Grenoble-1, France} 
\email{(Gaussier) herve.gaussier@ujf.grenoble.fr}
\address{(kim) Center for Geometry and its Applications 
and Mathematics Department, POSTECH, Pohang 790-784, The Republic of Korea}
\email{(Kim) kimkt@postech.ac.kr}
\thanks{Research of the first and the third named authors is 
supported in part by the grant 2011-0030044 (The SRC-GAIA) 
of the NRF of Korea.}%
%\subjclass[2010]{32A10}%

\begin{abstract}
The Kohn-Nirenberg domains are unbounded domains in 
$\CC^2$ upon which many outstanding questions are yet to be 
explored (cf., e.g.\ \cite{Fornaess_86}).  
The primary aim of this article is to demontrate that 
the Bergman and Caratheodory metrics of any Kohn-Nirenberg 
domain are positive and complete.
\end{abstract}
\maketitle

%% Section 1
\section{Introduction}

The following definition of the Kohn-Nirenberg domains is
due to Forn{\ae}ss \cite{Fornaess_86}.

\begin{definition}
\label{KN-domain}
A {\it Kohn-Nirenberg domain} (or, a {\it KN-domain} for the 
sake of brevity) is defined to be the set 
$$
\Omega = \{(z_1,z_2) \in \CC^2 \colon \re z_2 + P_{2k} (z_1) < 0 \},
$$ 
where $P_{2k}$ is a real-valued polynomial in $z_1$ and
$\bar z_1$ satisfying the following two conditions:
\begin{itemize}
\item[(\romannumeral1)] There exists an integer $k > 1$
such that $P_{2k}(rz) = r^{2k} P_{2k}(z)$ for any $r \in \RR$ 
and $z \in \CC$.
\smallskip

\item[(\romannumeral2)] $%\displaystyle 
\frac{\partial^2 P_{2k}}{\partial z_1 \partial \bar z_1}
\big|_{z_1} > 0$ for every $z_1 \neq 0$.
\end{itemize}
\end{definition}

The primary result of this paper is 

\begin{theorem} \label{main_thm}
The Bergman and Caratheodory metrics of any 
Kohn-Nirenberg domain are positive and complete.
\end{theorem}

Deferring further details to the next section, we clarify 
some terminology: (1) the {\it positivity} of a 
metric means that the length of any nonzero tangent vector is 
positive, and (2) the {\it completeness} of a metric means that 
the space equipped with its integrated distance is Cauchy-complete.
\smallskip

The study of KN-domains was initiated by the well-known work of Kohn 
and Nirenberg \cite{KN_73} which demonstrated the existence of 
pseudoconvex domains with a boundary point that does 
not admit, even locally, any holomorphic support functions
and consequently, the boundary cannot be convexifiable
there, via any change of holomorphic local coordinates.  
A representing example is the domain
$\Omega=\{(z_1,z_2) \in \CC^2 \colon \re z_2 + |z_1|^8 + 
\frac{15}7 |z_1|^2 ~\re z_1^6 < 0 \}$ with the
boundary point $(0,0) \in \partial\Omega$. 

Thereafter, many more questions have been posed on these 
domains, some of which still remain open. One remarkable 
result is by Forn{\ae}ss \cite{Fornaess_86}, which proves 
the sup-norm estimate for the $\bar\partial$-operator near 
the origin, a typical boundary point with no holomorphic 
support functions. 

It is still an open question whether every Kohn-Nirenberg 
domain is biholomorphic to a bounded domain. Another question 
along the same line was whether their Bergman metric is 
positive (`positive-definite', in some literatures) and 
complete. Notice that Theorem \ref{main_thm} answers this 
question affirmatively. \smallskip

Then, one recalls that \cite{KN_73} also introduces the domain
$$
W_\text{KN} := \{(z_1, z_2) \in \CC^2 \colon
\re z_2 + |z_1 z_2|^2 + |z_1|^8 + \frac{15}7 |z_1|^2 
~\re z_1^6 < 0 \}.
$$
This domain is special because it still does not allow any 
holomorphic local support function at the origin, despite
the fact that all of its boundary points except the origin 
are now strongly pseudoconvex. Notice however that this
modified defining function is not weighted-homogeneous 
any more. Similar domains with a degree 6 polynomial defining 
function of this type were discovered by Forn{\ae}ss 
\cite{Fornaess_77}; the domain
$$
W_\text{For} := \{(z_1, z_2) \in \CC^2 \colon \re z_1 + |z_1 z_2|^2 +
|z_1|^6 + t |z_1|^2 \re z_1^4 < 0\}
$$ 
for any $t$ satisfying $1<t<9/5$ enjoys the same property. 
Consequently, one naturally asks whether the conclusion of the 
theorem above countinues to hold for these domains. 
 
It turns out that the weighted homogeneity of the defining 
polynomial is not essential in the proof of Theorem 
\ref{main_thm}, and that a slight modification (cf.\ 
Proposition \ref{kn-org2}) yields

\begin{theorem} \label{main-2}
The Bergman and Caratheodory metrics of the domains
$W_\text{\rm KN}$ and $W_\text{\rm For}$ are positive and complete.
\end{theorem}

It seems to us that unbounded domains have emerged recently  and 
formed a territory in Several Complex Variables of high research 
interests (cf., e.g., \cite{HST}). In particular, it is not at all
clear whether the unbounded domains admit sufficiently many, 
independent, $L^2$ or $L^\infty$, holomorphic functions 
so that their Caratheodory and Bergman metrics can be seen to be 
positive, and furthermore, complete. 

For general complex manifolds, there are
well-known theorems based upon negativity of curvatures, e.g., 
\cite{Greene-Wu}.  But for domains in $\CC^n$, it is natural 
to ask whether their Bergman metric, for instance, can be seen to 
be positive and complete, directly from their defining functions.

One notable recent result in this direction is the following

\begin{theorem}[Chen-Kamimoto-Ohsawa \cite{CKO}]
If $\rho\colon \CC^n \to \RR$ is a nonnegative continuous 
plurisubharmonic function satisfying $\lim_{\|z\|\to\infty} 
\rho(z) = +\infty$, then the Bergman metric of the domain
$\Omega := \{(z, w) \in \CC^n \times \CC \colon 
\re w + \rho (z) < 0 \}$ 
in $\CC^{n+1}$ is positive-definite and complete.
\end{theorem}

Notice however that many Kohn-Nirenberg domains are not covered by 
this theorem.  Furthermore, the arguments and methods of
this paper are different from those of \cite{CKO}.

%% Section 2
\section{Basic observations for KN-domains}

\subsection{The Bergman kernel}
The {\it Bergman kernel} of a domain $W$ in $\CC^n$, denoted by
$K_W  (z, w)$, for $z, w \in W$, is constructed upon the 
Hilbert space $A^2(W)$ of $L^2$ holomorphic functions of 
$W$. (This space may happen to be trivial, however, since $W$ 
could be unbounded. In case $W$ is a complex manifold of
dimension $n$, $A^2(W)$ shall mean the Hilbert space of holomorphic 
$(n,0)$-forms, say $\alpha$, satisfying $\int_W \alpha \wedge 
\bar \alpha < +\infty$.) Then with the complete orthonormal system 
$\{\varphi_j \colon j=1,2,\ldots\}$ one has
$$
K_W (z,w) = \sum_{j=1}^\infty \varphi_j(z) \bar \varphi_j (w).
$$
(For the manifold case, 
$\sum \varphi_j(z) \wedge \bar \varphi_j (w)$.)
Then the {\it Bergman metric} is defined to be, if $K_W (z,z)>0$,
$$
\beta^W_z := \sum_{j,k=1}^n \frac{\partial^2\log K_W (z, z)}{\partial z_j \partial \bar z_k} \Big|_z \ dz_j \otimes d\bar z_k.
$$
While this $(1,1)$-tensor is known to define a positive-definite 
metric (a {\it positive} metric, in our terminology) if $W$ is a bounded
domain in $\CC^n$ by a theorem of Bergman himself, its positivity is
not at all clear even when $K_W (z, z) >0$ for every $z \in W$.  
Nevertheless, we shall call this tensor the Bergman {\it metric} 
following the usual convention.

In case the Bergman metric is positive, its real part defines a
Riemannian metric of $W$. Thus it generates the length of piecewise
$C^1$ curves and consequently the distance function, according to the 
usual routine of Riemannian geometry. We say therefore that the 
Bergman metric is complete if $W$, equipped with this distance, is 
Cauchy-complete as a metric space.
\medskip

Now let $\Omega := \{ (z_1, z_2) \in \CC^2 \colon 
\re z_2 + P_{2k} (z_1) < 0\}$ be a Kohn-Nirenberg domain 
in $\CC^2$ as in Definition \ref{KN-domain} above. Denote by
$K_\Omega$ the Bergman kernel function of $\Omega$. Then

\begin{proposition} 
\label{Bergman_kernel}
If $\Omega$ is a Kohn-Nirenberg domain, then it admits a nowhere
vanishing square-integrable holomorphic function. In particular, 
$K_\Omega (p, p) > 0$ for every $p \in \Omega$.
\end{proposition}

\noindent
\it Proof. \rm For this domain, there is a holomorphic 
peak function constructed by Bedford and Forn{\ae}ss  
\cite{Bedford-Fornaess} as follows:  
for sufficiently small a positive constant $\eta$, let 
\begin{equation} \label{BF_eta}
\Omega_\eta := 
\{(z_1, z_2)\in\CC^2 \colon \re z_2 
+ P_{2k} (z_1) < \eta \big(|z_2| + |z_1|^{2k} \big) \}.
\end{equation}
Then, by {\it Main Theorem} of Bedford and Forn{\ae}ss 
\cite{Bedford-Fornaess} p.\ 559, there exists a holomorphic function 
$\frak{f}\colon \Omega_\eta \to \CC$ satisfying:
\begin{enumerate}
\item There exists a constant $C>1$ such that
$$
\frac1C (|z_2| + |z_1|^{2k}) \le |\frak{f}(z)|
\le C (|z_2| + |z_1|^{2k}), ~\forall z=(z_1, z_2) \in \overline{\Omega}_\eta.
$$  
\item For a sufficiently large integer $N > 1$ there exists 
a branch of $\sqrt[N]{\frak{f}}$ such that 
$\arg \sqrt[N]{\mathfrak{f}} \in [-\pi/4, \pi/4]$.
\item $Q := \exp(-\sqrt[N]{\frak{f}})$ is a {\it holomorphic peak 
function} at $(0,0)$, i.e., $Q$ is holomorphic on $\Omega_\eta$, 
continuous on the closure $\overline{\Omega}_\eta$, $Q(0,0)=1$ 
and $|Q(z_1,z_2)| < 1$ for every
$(z_1,z_2) \in \overline{\Omega}_\eta \setminus \{(0,0)\}$.
\end{enumerate}
\medskip

Since this $Q(z)$ decays exponentially as $\|z\| \to \infty$ in $\Omega_\eta$, it is square-integrable on $\Omega_\eta$.  Since 
$\Omega \subset \Omega_\eta$, $Q$ is also square-integrable on 
$\Omega$. Since $Q$ does not vanish anywhere on $\Omega$, the 
proof follows immediately.
\hfill $\Box$
\medskip

Notice that the same conclusion holds for the Bergman kernel 
functions $K_{W_\textrm{KN}}$ and $K_{W_\textrm{For}}$ of the
domains $W_\textrm{KN}$ and $W_\textrm{For}$, respectively, since 
these domains are subdomains of KN-domains.

\subsection{The Caratheodory metric}

Recall the following classical concepts: for a complex 
manifold $M$, denote by $\cH(M, \Delta)$ the set of 
holomorphic functions from $M$ into the unit open disc $\Delta$ in 
$\CC$.  Let $p \in M$ and $v \in T_p M$. Then the 
{\it Caratheodory pseudo-metric} ({\it metric}, if positive) 
of $M$ is defined 
by
$$
F^C_M (p,v) 
= \sup \{ |df_p (v)| \colon f \in \cH(M, \Delta), f(p)=0 \}.
$$
This induces the {\it Caratheodory pseudo-distance} 
({\it distance}, if positive) 
$$
\rho^C_M (p, q) = \inf \int_0^1 F^C_M (\gamma(t), \gamma'(t)) dt,
$$
where the infimum is taken over all the piecewise $\cC^1$ curves
$\gamma\colon[0,1] \to M$ with $\gamma (0)=p$, $\gamma (1)=q$.
Of course, if we denote by $d^P_\Delta$ the Poincar\'e distance 
of the unit disc $\Delta = \{ z \in \CC \colon |z|<1\}$, then it
is well-known that $\rho^C_\Delta = d^P_\Delta$. 

\begin{remark}
Throughout this paper, we consider only the Caratheodory 
pseudo-distance introduced just now. 
But there is another equally well-known Caratheodory distance:
$$
d_\Omega^C (p,q) = \sup \{ d^P_\Delta (f(p), f(q)) \colon 
f \in \cH(M, \Delta) \}.
$$
We shall, however, deal with this pseudo-distance only almost 
at the end of this article, in Remark \ref{fin_rem}. We point 
out on the other hand, in most other literatures, our Caratheodory 
pseudo-distance is usually called the {\it integrated} Caratheodory 
pseudo-distance.
\end{remark}

Now we present

\begin{proposition}
If $\Omega$ is a Kohn-Nirenberg domain, then
$F^C_\Omega (p,v) > 0$ for every $(p,v) \in \Omega \times 
\big(\CC^2 \setminus \{(0,0)\}\big)$. 
\end{proposition}

\noindent\it Proof. \rm
For $p = (p_1, p_2)$ and $v=(v_1,v_2)$ take $g(z_1, z_2) := 
Q(z_1, z_2) \cdot \big(\bar v_1 (z_1-p_1)
+ \bar v_2 (z_2-p_2) \big)$.  Then $g$ is a bounded holomorphic 
function of $\Omega$, since $Q$ decays exponentially at infinity in 
$\Omega$ and is continuous in the closure. Moreover 
$g(p) = 0$ and $|dg_p (v)| = Q(p)\|v\|^2 > 0$.  Hence the proof is
complete.
\hfill $\Box$ 

\subsection{The Hahn-Lu comparison theorem}

The following theorem is very mild a modification of the comparison
theorem by Hahn \cite{Hahn1, Hahn_77, Hahn2}, and Lu \cite{Lu} which
compares the Caratheodory metric and the 
Bergman metric:

\begin{theorem}[Hahn-Lu comparison theorem] 
If $M$ is a complex manifold such that 
\begin{itemize}
\item[(\romannumeral 1)] its Caratheodory metric $F^C_M$ is
positive, and 
\item[(\romannumeral 2)] its Bergman kernel $K_M$ satisfies
$K_M (p, p) \neq 0$ for every $p \in M$,
\end{itemize}
then its Bergman metric $\beta^M_p (v, w)$ satisfies 
the inequality
$$
\big(F^C_M (p, v)\big)^2 \le \beta^M_p (v, v), 
$$
for any $p \in M$ and $v \in T_p M$. In particular, this implies
that the Bergman metric is positive.
\label{HL}
\end{theorem}

\noindent\it Proof. \rm We shall only prove it for the 
case when $M=\Omega$
is a domain in $\CC^n$, fitting to the purpose of this article;  
the manifold case uses essentially the same arguments except some
simplistic adjustments. 
  
Start with the following quantities developed
by Bergman \cite{Bergman_35}: 
\begin{eqnarray*}
{\frak B}_0 (p) 
& = & \sup\Big\{|\psi (p)|^2 \colon \psi 
\hbox{ holomorphic}, \int_\Omega |\psi|^2 \le 1 \Big\}
\\
{\frak B}_1 (p, v) 
& = & \sup\Big\{ \Big|\partial_v \varphi |_{p}\Big|^2 
\colon \varphi \textrm{ holomorphic}, \varphi(p)=0,
\int_\Omega |\varphi|^2 \le 1 \Big\},
\end{eqnarray*}
where $%\displaystyle 
\partial_v \varphi|_p = \sum_{j=1}^n v_j 
\frac{\partial \varphi}{\partial z_j}\Big|_{p}$. 
These concepts are significant because 
${\frak B}_1 (p, v) = {\frak B}_0 (p) \cdot
\beta^\Omega_p (v, v)$, when ${\frak B}_0 (p)>0$.

Following \cite{Hahn2}, consider an $L^2$-holomorphic 
function $\hat\psi$  on $\Omega$ with 
$\|\hat\psi\|_{L^2(\Omega)} \le 1$ satisfying
$|\hat\psi (p)|^2 = \frak{B}_0(p)$.
Then Montel's theorem on normal families implies the existence of
$\eta \in \cH(\Omega, \Delta)$ on $\Omega$ with 
$\eta(p)=0$ and $\big|\partial_v \eta|_{p}\big|=
|d\eta_{p} (v)| = F_\Omega^C (p, v)$, the 
Caratheodory length of $v$ at $p$.  Since
$|\eta\hat\psi| \le |\hat\psi|$, one obtains
$\big|\partial_v  (\eta\hat\psi)|_{p}\big|^2 
\le \frak{B}_1 (p,v)$.
Since the function $\eta\hat\psi$ is holomorphic on 
$\Omega$ with $\|\eta\hat\psi\|_{L^2(\Omega)} \le \|\hat\psi\|_{L^2(\Omega)} \le 1$ and vanishes at $p$, 
one immediately obtains that
$\big|\partial_v  (\eta\hat\psi)|_{p}\big|
= \big|\partial_v  \eta|_{p}\big|~ |\hat\psi(p)|$
and more importantly, that
$\frak{B}_1 (p,v) \ge F^C_\Omega (p, v)^2 \ 
\frak{B}_0(p).$
Hence the comparison
$$
\beta^{\Omega}_p (v, v) = \frac{\frak{B}_1 (p,v)}{\frak{B}_0(p)} \ge 
F^C_\Omega (p, v)^2
$$
follows as desired.
\hfill $\Box$
\medskip

The original statement required positivity of both 
metrics. But the proof above (with no changes from the 
arguments by Hahn \cite{Hahn2}) clearly shows that not all those 
assumptions are necessary. In fact these thoughts yield:

\begin{proposition}
If $W$ is a Kohn-Nirenberg domain in the sense of Definition 
\ref{KN-domain} or one of the domains $W_\text{\rm KN}$ and 
$W_\text{\rm For}$, then their Bergman and Caratheodory  
metrics are positive. Moreover, if we denote by $\beta^W_p (v, w)$ 
the Bergman metric at $p$, then
$$
\beta^W_p (v, v) \ge F^C_W (p,v)^2
$$
for every $(p ,v) \in W \times \CC^2$.
\end{proposition}

\section{Construction of peak functions and Completeness}

We establish, in this section, the completeness of the 
Bergman and Caratheodory metrics of the Kohn-Nirenberg domains. 
Thanks to the comparison theorem 
of Hahn-Lu, we are only to show that any KN-domain (as well as
$W_\text{\rm KN}$ and $W_\text{\rm For}$) equipped
with its Caratheodory distance is Cauchy-complete.
\medskip

Towards this goal, the following statement plays a crucial role.

\begin{theorem} \label{peaks}
Every boundary point of the Kohn-Nirenberg domain $\Omega$
admits a holomorphic peak function. More precisely, every 
boundary point $p \in \partial\Omega$ admits a holomorphic 
function $f_p\colon \Omega \to \Delta$ satisfying:
\begin{itemize}
\item[(1)] $\lim_{z \to p} f_p(z) = 1$.
\item[(2)] For any $r>0$ there exists $s>0$ such that 
$|f_p (w)|<1-s$ whenever $w \in \Omega \setminus B(p,r)$. 
\end{itemize}
\end{theorem}

\noindent\it Proof. \rm
As before, our Kohn-Nirenberg domain $\Omega$ is defined by
the inequality $\re z_2 + P_{2k} (z_1) < 0$.
We construct a holomorphic peak function, say 
$f_p \in \cH(\Omega, \Delta)$, at every boundary point 
$p \in \partial\Omega$.  

If $p=(0,0)$, this is already done in 
\cite{Bedford-Fornaess}; the function $Q$ in the proof of
Proposition \ref{Bergman_kernel} suffices. Then the 
translation along the $\im z_2$ direction yields holomorphic 
peak functions at the points $p = (0, ib) \in \partial\Omega$ 
for every $b \in \RR$.
\smallskip

If $p \in \partial\Omega \setminus \{(0, ib) \colon b \in \RR\}$, 
$\partial\Omega$ is strongly pseudoconvex at $p$. 
Then there is of course a holomorphic local-peak function,
i.e., a holomorphic function $\mathfrak{g}_p\colon B(p,r) \to 
\CC$ for some $r>0$ such that $\mathfrak{g}_p(p)=1$ and 
$|\mathfrak{g}_p (z)| < 1$ for every 
$z \in B(p,r) \cap \overline{\Omega} \setminus \{p\}$.
However, since the domain is unbounded, care has to be taken
in extending $\mathfrak{g}_p$ to a global peak function.

Take constants $0<r_1<r_2<r$ and then choose $\varepsilon >0$
so that the domain
$
\Omega^\varepsilon = \{ z\in \CC^2 \colon
\re z_2 + P_{2k} (z_1) < \varepsilon \}
$
is strongly pseudoconvex at every
$q \in \partial \Omega^\varepsilon \cap B(p,r_2) $ and that
$$
\{z \in B(p,r) \colon \mathfrak{g}_p(z)=1\} \cap 
\big(B(p,r_2) \setminus \overline B(p,r_1) \big)
\subset \CC^2 \setminus \overline{\Omega^\varepsilon}.
$$
Then take the $\cC^\infty$ cut-off function
$\chi \colon \CC^2 \to [0,1]$ satisfying 
$\chi \equiv 1$ on $B(p,r_1)$ and $\supp \chi \subset B(p,r_2)$.

Now we exploit the $\bar\partial$-problem on 
$\Omega^\varepsilon$ with an $L^2$-estimate with an arbitrary 
plurisubharmonic weight (\cite{Hor_73}, Theorem 4.4.2). 
First let
$$
\alpha_p (z) = \begin{cases} 
\bar\partial\Big(\frac{\chi(z)}{(1-\mathfrak{g}_p(z))Q(z)} \Big)
& \textrm{if } z \in \Omega^\varepsilon \cap (B(p,r_2) \setminus 
\overline B(p,r_1)) 
\\
0 & \textrm{if } z \in \big(\Omega^\varepsilon \setminus 
\overline B(p,r_2)\big) 
\cup B(p,r_1)).
\end{cases} 
$$
(We acknowledge at this point that the set-up of this
differential form has been influenced by the methods 
of Forn{\ae}ss and McNeal \cite{FM}.) 

Then consider the function $u_p\colon\Omega^\varepsilon\to\CC$ that solves 
the equation $\bar\partial u_p = \alpha_p$ with the estimate
(with the zero plurisubharmonic weight)
$$
\int_{\Omega^\varepsilon} \frac{|u_p(z)|^2}{(1+\|z\|^2)^2} 
d\mu(z)
\le
\int_{\Omega^\varepsilon} |\alpha_p(z)|^2 d\mu(z).
$$
Notice that this $\alpha_p$ is a bounded-valued
$\bar\partial$-closed smooth $(0,1)$-form with bounded 
support inside $\Omega^\varepsilon$. Consequently the 
right-hand side of the preceding inequality is bounded, say,
by a positive constant $A$.
\smallskip

We would like to obtain a pointwise estimate for $|u_p(z)|$ 
at every $z\in\Omega$ with $\|z\|\gg 1$. 
Take a positive constant $R_p$ sufficiently large that 
$B(p,r_2) \subset B(0,R_p-1)$. Since the defining function is
a polynomial, one sees that,
for any $\xi \in \Omega \setminus B(0, R_p)$, there exists a uniform 
$c_\varepsilon>0$ independent of $\xi$ such that 
$B(\xi, c_\varepsilon\|\xi\|^{-2k})\subset \Omega^\varepsilon$ 
and that the support of the function $\chi$ has no 
intersection with $B(\xi, c_\varepsilon\|\xi\|^{-2k})$. 
In particular, $u_p$ is holomorphic on 
$B(\xi, c_\varepsilon\|\xi\|^{-2k})$.

Thus, for any $\xi=(\xi_1, \xi_2) \in \Omega \setminus B(0,R_p)$,
\begin{eqnarray*}
A & \ge & \int_{\Omega^\varepsilon} |\alpha_p(z)|^2 d\mu(z) 
\\
& \ge & 
\int_{\Omega^\varepsilon} \frac{|u_p(z)|^2}{(1+\|z\|^2)^2} d\mu(z)
\\
& \ge & 
\int_{B(\xi, c_\varepsilon\|\xi\|^{-2k})} 
\frac{|u_p(z)|^2}{(1+\|z\|^2)^2} d\mu(z)
\\
& \ge & 
\frac{1}{4\|\xi\|^4}\int_{B(\xi, c_\varepsilon\|\xi\|^{-2k})} |u_p(z)|^2 d\mu(z).
\end{eqnarray*}
Since $u_p$ is holomorphic on 
$B(\xi, c_\varepsilon\|\xi\|^{-2k})$, the sub mean-value
inequality implies that $|u_p(\xi)|^2 \le A_p\|\xi\|^{8k+4}$,
for some positive constant $A_p$ independent of $\xi$.
In particular, $|u_p|$ grows at most polynomially (of
degree $8k+4$) at infinity.

Notice that $u_p$ is smooth in $\Omega^\varepsilon$ and hence
smooth on the closure of $\Omega$.  Since the 
Bedford-Forn{\ae}ss peak function $Q$ enjoys 
the exponential decay estimate at infinity, we may take
sufficiently small a constant $c_p>0$ so that
$$
c_p |Q(z) u_p (z)| < \frac12, \quad 
\forall z \in \Omega.
$$
In particular, we have
\begin{equation} \label{at_infty}
\re \big(c_p Q(z) u_p (z) - 1\big) < -\frac12, 
\quad \forall z \in \Omega.
\end{equation}
Since $u_p$ satisfies the equation 
$\bar\partial u_p = \alpha_p$, it follows that 
the function
$$
\frac{\chi(z)}{1-\mathfrak{g}_p(z)} -  Q(z) u_p(z)
$$ 
is holomorphic on $\Omega$.
\smallskip

Altogether, if we define $f_p\colon \Omega \to \CC$ with the positive 
constant $c_p$ by 
$$
f_p (z) = \exp \left(\frac{\mathfrak{g}_p(z)-1}{c_p \chi (z) - 
(1-\mathfrak{g}_p(z)) [c_p Q(z) u_p (z) - 1]}\right),
$$
then $f_p \in \cA(\Omega)$.  Moreover, for the positive constant 
$c_p$, one easily obtains by \eqref{at_infty} that
$$
\re \left(\frac{\mathfrak{g}_p(z)-1}{c_p \chi (z) - 
(1-\mathfrak{g}_p(z)) [c_p Q(z) u_p (z) - 1]}\right) < 0
$$
for every $z \in \Omega$.

Finally, $\lim_{\Omega \ni z\to p} f_p (z) = e^0 = 1$. 
Notice that $p$ is the only boundary point
that has this property for $f_p$. 
The other condition is also easily checked since this last
estimate receives a definite negative upper bound if $z$
is at a positive distance away from $p$. 
Therefore, $f_p$ is the desired global peak function at $p$
for $\Omega$. \hfill $\Box$
\medskip

\begin{remark} \label{KN-For}
\rm
Notice that the homogeneity of the polynomial $P_{2k}$ is 
not essential in extending the holomorphic local-peak function
$\frak{g}_p$ to the global peak function $f_p$ when 
$p \in \partial\Omega$ is a strongly pseudoconvex point. 
The same extension procedure works for any domain  
$W_S := \{(z_1, z_2) \in \CC^2 \colon 
\re z_1 + S(z_1, z_2) + P_{2k} (z_1) < 0 \}$ 
defined by the polynomial defining function, whenever $P_{2k}$ 
is as in Definition \ref{KN-domain}, as long as $S$ satisfies the 
following technical conditions:
\begin{itemize}
\item[(\romannumeral1)] $S$ is a nonnegative real-valued 
polynomial in $z_1, \bar z_1, z_2$, and $\bar z_2$ of degree 
strictly less than $2k$ with $S(0)=0$; and
\item[(\romannumeral2)]  $W_S$ is strongly pseudoconvex at every boundary 
point possibly except the origin. 
\end{itemize}
Notice that the domains $W_\text{KN}$ and
$W_\text{For}$ belong to such a collection of domains. Thus the 
conclusion of the theorem above holds in particular for the
domains $W_\text{KN}$ and $W_\text{For}$.
\end{remark}

Now we prove

\begin{theorem} \label{Carath-complete}
Every Kohn-Nirenberg domain defined in Definition \ref{KN-domain}, 
equipped with its Caratheodory distance, is a 
complete metric space.
\end{theorem}

\noindent\it Proof. \rm
Let $(q_n)_n$ be a Cauchy sequence in the Kohn-Nirenberg domain
$\Omega$ with respect to the Caratheodory distance 
$\rho_\Omega^C$. We now pose:
\medskip

\bf Claim. \it There exists a compact subset $X$ of 
$\CC^2$ such that $q_n \in X \subset \Omega$ for every $n=1,2,\ldots$.
\rm
\medskip

Note that this implies that $(q_n)_n$ has a subsequence 
$(q_{n_k})_k$ convergent to $\hat q \in X$ with respect to the 
Euclidean distance. However, since the 
Caratheodory distance $\rho_\Omega^C$ is continuous, 
$(q_{n_k})_k$ converges to $\hat q$ with respect to $\rho_\Omega^C$, and 
hence $(q_n)_n$ converges to $\hat q$ with respect to $\rho_\Omega^C$
as well. Therefore, it suffices to establish this claim.

Theorem \ref{peaks} implies that this sequence has no subsequence 
that approaches a boundary point of $\Omega$ arbitrarily closely 
in the Euclidean distance, due to the distance-decreasing property 
of the Caratheodory distance. Hence it remains 
to show that no subsequence of $(q_n)_n$ can diverge indefinitely 
far away from the origin with respect to the Euclidean distance. 

Assume the contrary.  Then choosing a subsequence we may
assume without loss of generality that
$\lim_{n\to\infty} \|q_n\| = \infty$.  Let $q_n := (a_n, b_n)$
for every $n$. Notice that, for every $n=1,2,\ldots$, there exists a 
unique positive number $t_n$ satifying the equation
$$
t_n^2 |a_n|^2 + t_n^{4k} |b_n|^2 = 1.
$$
Define the map $\varphi_n\colon \CC^2 \to \CC^2$ by 
$\varphi_n (z_1, z_2) = (t_n z_1, t_n^{2k} z_2)$ for each $n$.
Then $\varphi_n \in \aut\Omega$. 

Observe that $\lim_{n\to\infty} \varphi_n (q_1) = (0,0)$.
If we denote by $S = \{(z_1,z_2) \in \CC^2 \colon 
|z_1|^2 + |z_2|^2 = 1 \}$, then $\varphi_n (q_n) \in 
S \cap \Omega$, for every $n$. Therefore,
$$
\rho_\Omega^C (q_1, q_n) = \rho_\Omega^C (\varphi_n (q_1), 
\varphi_n (q_n)) \ge 
\rho_\Delta^C \big(Q(\varphi_n (q_1), Q(\varphi_n (q_n)\big), 
$$
and $\|\varphi_n (q_1) - \varphi_n (q_n)\| > \frac12$ for $n\gg 1$.
Then Proposition \ref{peaks} yields a constant $s$ with $0<s<1$ 
independent of $n$ with $|Q(\varphi_n (q_n))|<1-s$. 
Therefore one sees that
$$
\rho_\Delta^C \big(Q(\varphi_n (q_1), Q(\varphi_n (q_n)\big)
\ge
\inf\limits_{t\in\RR} \rho_\Delta^C \big(Q(\varphi_n (q_1)), 
(1-s)e^{it}\big).
$$
But, since $\lim_{n\to\infty} \varphi_n (q_1) 
= (0,0)$, the preceding inequality yields that
$$
\lim_{n\to\infty} \rho_\Omega^C (q_1, q_n) = +\infty.
$$ 
So $(q_n)_n$ fails to be a bounded sequence with respect to
the Caratheodory distance, and consequently
cannot be a Cauchy sequence with respect to the Caratheodory 
distance $\rho^C_\Omega$. 
This contradicts the original hypothesis that $(q_n)_n$ was
$\rho^C_\Omega$-Cauchy. Thus the proof is complete. 
\hfill $\Box$
\medskip

Furthermore, we obtain

\begin{proposition} \label{kn-org2}
Let $\Omega := \{(z_1,z_2) \in \CC^2 \colon \re z_2 
+ P_{2k}(z_1)<0 \}$ be a Kohn-Nirenberg domain. If $S(z_1, z_2)$ 
is a nonnegative polynomial of degree strictly less than $2k$, 
with $S(0)=0$, such that the domain 
$$
W_S := \{(z_1,z_2) \in \CC^2 \colon 
\re z_2 + S(z_1, z_2) + P_{2k} (z_1) < 0\}
$$
is strictly pseudoconvex everywhere except at the origin $(0,0)$, 
then the domain $W_S$ equipped with its Caratheodory 
distance is Cauchy-complete.
\end{proposition}

\noindent\it Proof. \rm 
Since $W_S$ is contained in the Kohn-Nirenberg domain 
$\Omega$, $\rho^C_{W_S} \ge \rho^C_\Omega$ and consequently
$\rho^C_{W_S}$ is positive. By Remark \ref{KN-For}, $W_S$ admits
a holomorphic peak function at every strongly pseudoconvex point. 
Since the origin is a common boundary point to $W_S$ and $\Omega$, the 
peak function for $\Omega$ at the origin is also a peak function of 
$W_S$ at the origin. Thus no Cauchy sequence of $W_S$ with respect
to its Caratheodory distance $\rho_{W_S}^C$ can 
accumulate at a boundary point of $W_S$.  Finally, since 
the Caratheodory distance $\rho^C_{W_S}$ is larger than 
or equal to $\rho^C_\Omega$, no Cauchy sequence 
with respect to $\rho^C_{W_S}$ can diverge
indefinitely far away from the origin.  Altogether, it follows
that every Cauchy sequence of $W_S$ is bounded, and bounded away 
from the boundary. Now the proof follows by the continuity of the
Caratheodory distance. 
\hfill $\Box$
\bigskip

Altogether the proofs of Theorems \ref{main_thm} and \ref{main-2} 
now follow by the Hahn-Lu comparison theorem (Theorem \ref{HL}).  
Of course the domains described in this proposition also have their 
Bergman metrics positive and complete.

\section{Higher Dimensions}

The method of this paper up to this point is not restricted 
to dimension two, except at the places where 
the Bedford-Forn{\ae}ss peak functions were exploited. 
But that is also not a strong restriction. In fact the following 
existence theorem of peak functions at the origin of certain 
domains in higher dimensional cases were established by Noell 
\cite{Noell}:
\medskip

Let $U$ be a domain in $\CC^n$. 
Denote by $z = (z', z_n)$ the complex variable(s) for $\CC^n$
with $z'=(z_1, \ldots, z_{n-1})$, and by $\cA(U)$ the set of 
continuous functions on $\overline U$, holomorphic on the domain $U$. 

\begin{theorem}[Noell]
\label{highdim_peak}
Suppose that $P_{2k}$ is a homogeneous plurisubharmonic 
polynomial of degree $2k$ on $\CC^{n-1}$, and assume that $P_{2k}$ 
is not harmonic along any complex line through the origin of 
$\CC^{n-1}$. For a sufficiently small positive constant $\eta$, let
$$
\Omega_\eta:=\{(z', z_n) \in \CC^n \colon \re z_n+ P_{2k}(z')
<\eta|z_n|+\eta\|z'\|^{2k}\}.
$$
Then, there exists a function $\mathfrak{f}\in\cA(\Omega_\eta)$ 
with the following properties:
\begin{enumerate}
\item For some constant $C>1$,
$$
\frac1C (|z_n| + \|z'\|^{2k}) \leq |\frak{f}(z)|
\leq C(|z_n| + \|z'\|^{2k}), 
$$ 
for all $z=(z', z_n) \in \overline{\Omega}_\eta$.
\item For a sufficiently large integer $N > 1$ there exists 
a branch of $\sqrt[N]{\frak{f}}$ such that 
$\arg \sqrt[N]{\mathfrak{f}} \in [-\pi/4, \pi/4]$.
\item For all $i_1, \ldots, i_n \geq 0$ there exist constants 
$C_{i_1, \ldots, i_n}$ and $N_{i_1, \ldots, i_n}$ such that
$$
\Big|\frac{\partial^{i_1+ \ldots +i_n}\frak{f}}
{\partial z_1^{i_1} \ldots \partial z_n^{i_n}}\Big| 
\leq C_{i_1, \ldots, i_n} \|z\|^{-N_{i_1, \ldots, i_n}}
\quad \textrm{ when } \|z\| \leq 1.
$$
\item $\exp(-\sqrt[N]{\frak{f}})$ is a peak function at $0$ 
for $\cA(\Omega_\eta)$.
\item $\exp(-1/\sqrt[N]{\frak{f}}) \in C^\infty 
(\overline{\Omega}_\eta)$ and $\frak{f}$ is a separating 
function at the origin for $\cA(\Omega_\eta)$.
\end{enumerate}
\medskip
\end{theorem}

Thus we immediately obtain

\begin{theorem}
The Bergman metric and the Caratheodory metric of 
any domain $\Omega$ in the statement of the preceding theorem are 
positive and complete.
\end{theorem}

Of course, all other theorems in this article also generalize to 
these domains.

\begin{remark} \rm \label{fin_rem}
Another well-known {\it Caratheodory pseudo-distance} is defined by
$$
d^C_M (p,q) = \sup \{d^P_\Delta (f(p), f(q)) \mid 
f \colon M \to \Delta \text{ holomorphic} \},
$$
where $d^P_\Delta$ is the Poincar\'e distance of the unit disc
$\Delta$.  This Caratheodory distance $d^C_M$ is in general
smaller than the (integrated) Caratheodory distance $\rho^C_\Omega$. 
Hence one might like to ask whether KN-domains as well as
the higher dimensional domains just mentioned here, when equipped 
with the Caratheodory distance $d^C$, are also complete.  
They indeed are. 

Using peak functions one can show that the Cauchy sequences 
cannot approach the boundary. The sequence cannot diverge 
infinitely far from the origin by the proof-arguments of 
Theorem \ref{Carath-complete}.  The only remaining point to show now
should be the positivity, i.e., $d_\Omega^C (p,q) >0$ if $p\neq q$. 
But this was established earlier by Yu \cite{Yu_95}.
\end{remark}

%% References

\end{document}